\documentclass[a4paper,12pt]{article} 
\usepackage{amssymb}
\usepackage{amscd}
\usepackage{amsmath} 
\usepackage[dvipdfmx]{graphicx}
\usepackage{latexsym} 
\usepackage{enumerate}
\usepackage[usenames]{color}
\usepackage{lscape}

\usepackage{theorem}
\theoremstyle{plain}
\theorembodyfont{\itshape}
\newtheorem{theorem}{Theorem}[section]
\newtheorem{proposition}[theorem]{Proposition}
\newtheorem{lemma}[theorem]{Lemma}

\theorembodyfont{\rmfamily}

\newtheorem{remark}[theorem]{Remark}

\newtheorem{recipe}[theorem]{Recipe} 

\def\ri{\mathrm{i}}

\def\rM{\mathrm{M}}
\def\diag{\mathrm{diag}}

\def\rRe{\mathrm{Re}} 
\def\C{\mathbb{C}}

\def\Q{\mathbb{Q}}
\def\R{\mathbb{R}}
\def\Z{\mathbb{Z}}

\def\cM{\mathcal{M}}
\def\cS{\mathcal{S}}
\def\ve{\varepsilon}
\def\vG{\varGamma}
\def\vD{\varDelta}
\def\vL{\varLambda}
\def\ba{\mbox{\boldmath $a$}}
\def\sba{\mbox{\boldmath ${\scriptstyle a}$}}

\def\bc{\mbox{\boldmath $c$}}

\def\be{\mbox{\boldmath $e$}}
\def\bk{\mbox{\boldmath $k$}}

\def\bp{\mbox{\boldmath $p$}}
\def\sbp{\mbox{\boldmath ${\ss p}$}}
\def\bq{\mbox{\boldmath $q$}}
\def\br{\mbox{\boldmath $r$}}

\def\bu{\mbox{\boldmath $u$}}
\def\bv{\mbox{\boldmath $v$}}
\def\bw{\mbox{\boldmath $w$}}

\def\0{\mbox{\boldmath $0$}}
\def\1{\mbox{\boldmath $1$}}
\def\2{\mbox{\boldmath $2$}}
\def\3{\mbox{\boldmath $3$}}

\def\carl{\circlearrowleft}
\def\hgf{{}_3f_2}
\def\hgF{{}_3F_2}

\def\hgh{{}_3h_2}

\def\bhgh{{}_3\mbox{\boldmath $h$}_2}
\def\pFq{{}_pF_q}
\def\ds{\displaystyle}
\def\ss{\scriptstyle} 
\def\ts{\textstyle}
\title{\bf Three-Term Relations for $\mbox{\boldmath ${}_3F_2(1)$}$}   
\author{Akihito Ebisu\thanks{Faculty of Information and Computer Science, 
Chiba Institute of Technology 2-1-1, Shibazono, Narashino, Chiba, 275-0023, Japan. 
{\tt akihito.ebisu@p.chibakoudai.jp}} \ and 
Katsunori Iwasaki\thanks{Department of Mathematics, 
Hokkaido University, Kita 10, Nishi 8, Kita-ku, Sapporo 060-0810 Japan. 
{\tt iwasaki@math.sci.hokudai.ac.jp} (Corresponding author)}} 
\date{September 4, 2017} 
\begin{document}
\maketitle
\begin{abstract} 
For the hypergeometric function of unit argument $\hgF(1)$ we prove  
the existence and uniqueness of three-term relations with arbitrary 
integer shifts. 
We show that not only the original $\hgF(1)$ function but also other 
five functions related to it satisfy one and the same three-term relation. 
This fact is referred to as simultaneousness. 
The uniqueness and simultaneousness provide three-term relations 
with a group symmetry of order $72$.       
\end{abstract} 
\section{Introduction} \label{sec:intro}
If $p \le q+1$ then the hypergeometric function $\pFq$ admits $(q+2)$-term 
contiguous relations (see Rainville \cite{Rainville}). 
Under certain conditions, they may reduce to ones with only a smaller 
number of terms.  
In the case of $\hgF$ the general contiguous relations are 
four-term ones, while Kummer observed that it was possible to obtain 
three-term contiguous relations for $\hgF$ when the argument was $1$, 
that is, for $\hgF(1)$; see Andrews {\sl et al.} \cite[\S 3.7]{AAR}.    
Bailey \cite{Bailey1} gave a procedure to produce those relations 
using differential equations and Wilson \cite{Wilson} gave a 
simpler method.    
\par
A contiguous function in its strict sense is a function obtained 
from the original $\pFq$ by altering one of the parameters by $\pm1$. 
A three-term contiguous relation is a linear relation between  
the original $\pFq$ and two other functions contiguous to it. 
In the case of $\hgF(1)$ there are a total of twelve contiguous 
relations, excluding the ones obtained by permuting numerator 
or denominator parameters, a complete list of which can be 
found in \cite[formulas (13)--(24)]{Wilson}. 
\par
We can generalize the concept of contiguity by calling a 
function contiguous if it is obtained from 
the original function by shifting the parameters by 
{\sl an arbitrary integer vector}.  
We can then speak of a {\sl general three-term relation}  
as a linear relation among three functions contiguous 
to each other in the generalized sense. 
We may safely say that the study of three-term contiguous 
relations for $\hgF(1)$ in the narrow sense is finished by 
the works of Bailey \cite{Bailey2} and Wilson \cite{Wilson}.  
This is not the case with the general three-term relations; 
there is still something very basic that should be added to 
the subject. 
We establish such results (Theorem \ref{thm:3tr}) regarding   
the existence, uniqueness and simultaneousness of general 
three-term relations for $\hgF(1)$, where simultaneousness 
means that one three-term relation is commonly satisfied 
by certain six functions associated with the original 
$\hgF(1)$ function. 
As a corollary to the uniqueness and simultaneousness we 
are also able to obtain a group symmetry of order $72$ 
on three-term relations. 
Underlying the existence and uniqueness is the linear independence 
over the rational function field of two $\hgF(1)$ functions 
whose variables differ by a nonzero integer vector. 
\par
For the Gauss hypergeometric function ${}_2F_1(\ba; z) := 
{}_2F_1(a, b; c; z)$ with parameters $\ba := (a,b;c) \in \C^3$, 
Vidunas \cite{Vidunas} considered three-term relations representing 
${}_2F_1(\ba + \bk; z)$ for $\bk \in \Z^3$ in terms of 
${}_2F_1(\ba; z)$ and ${}_2F_1(\ba + \be_1; z)$, where $\be_1 := (1,0;0)$. 
He showed the existence and uniqueness of three-term relations of 
this form \cite[Theorem 1.1]{Vidunas} and obtained simultaneousness 
for them \cite[formulas (19)--(23)]{Vidunas}. 
There is a similar approach that carries over when the particular 
shift vector $\be_1$ is replaced by an arbitrary nonzero integer vector. 
To discuss similar issues for $\hgF(1)$, however, we shall develop a quite  
different method that works for every integer shift.   
As for ${}_2F_1$, Ebisu \cite{Ebisu1,Ebisu2} gave a useful formula for 
three-term relations expressing ${}_2F_1(\ba + \bk; z)$ in terms of 
${}_2F_1(\ba; z)$ and ${}_2F_1(\ba + \1; z)$ with $\1 := (1,1;1)$, 
derived symmetry on them from simultaneousness and moreover applied 
these results to special values.         
\par
Recall that the hypergeometric series $\hgF(\ba; z)$ 
is a power series of $z$ defined by     
\begin{equation} \label{eqn:hgF(z)}
\hgF(\ba; z) := \sum_{k=0}^{\infty} 
\dfrac{(a_0, k) \, (a_1, k) \, (a_2, k) }{(1, k) \, (a_3, k) \, (a_4, k)} 
\, z^k,  
\qquad (a, k) := \frac{\vG(a+k)}{\vG(a)},   
\end{equation}
with complex parameters $\ba = (a_0,a_1,a_2;a_3,a_4) \in \C^5$, 
which are often denoted by 
\[
\ba := 
\begin{pmatrix}
a_0, & a_1, & a_2 \\
     & a_3, & a_4 
\end{pmatrix}
= 
\begin{pmatrix}
a_0, & a_1, & a_2 \\
     & b_1, & b_2 
\end{pmatrix}.   
\]
Throughout we will freely switch between $(a_3, a_4)$ and 
$(b_1, b_2)$ according to the situation. 
\par
In this article it is more convenient to work with a 
renormalized version of the series \eqref{eqn:hgF(z)}:     
\begin{equation} \label{eqn:hgf(z)}
\hgf(\ba; z) 
:= 
\sum_{k=0}^{\infty} 
\dfrac{\vG(a_0 +k) \, \vG(a_1 +k) \, \vG(a_2 + k) }{\vG(1+k) \, \vG(a_3 +k) \, \vG(a_4 +k)} 
\, z^k 
= \frac{\vG(a_0) \, \vG(a_1) \, \vG(a_2)}{\vG(a_3) \, \vG(a_4)} \, \hgF(\ba; z).        
\end{equation}
If $\ba$ satisfies $a_i \not \in \Z_{\le 0}$, $i = 0,1,2,3,4$, then 
\eqref{eqn:hgf(z)} is a well-defined power series in $z$ with 
non-vanishing leading coefficient. 
To ensure this condition for all $\ba + \bp$ with $\bp \in \Z^5$ we 
assume    
\begin{equation} \label{eqn:wd}
a_0, \, a_1, \, a_2, \, a_3, \, a_4 \not \in \Z \qquad \mbox{for} 
\quad \ba = (a_0,a_1,a_2;a_3,a_4) \in \C^5.   
\end{equation}
\par
It is well known that $\hgf(\ba; z)$ converges in $|z| < 1$ and 
solves a Fuchsian differential equation    
\begin{equation} \label{eqn:ode}
L(\ba) \, y := \{ \, \theta (\theta+b_1-1)(\theta+b_2-1) - 
z \, (\theta +a_0)(\theta +a_1)(\theta +a_2) \, \} \, y = 0,  
\quad \theta := z \ts \frac{d}{d z},    
\end{equation}
the third-order hypergeometric equation, whose Riemann scheme is 
given by 
\begin{equation} \label{eqn:rs}
\left\{
\begin{matrix}
z = 0 & z = 1  & z = \infty \\
1-b_0 & 0      & a_0        \\
1-b_1 & 1      & a_1        \\
1-b_2 & s(\ba) & a_2 
\end{matrix}
\right\}, \qquad 
s(\ba) := b_1 + b_2 - a_0 - a_1 - a_2,   
\end{equation}
where $b_0 := 1$ by convention and $s(\ba)$ is referred to as the 
{\sl Saalsch\"utzian index} of $\ba$. 
\par
\begin{table}[t]
\begin{align*}
\sigma_0^{(0)}(\ba) &:= \ba =  
\begin{pmatrix}
a_0, & a_1, & a_2 \\
     & b_1, & b_2 
\end{pmatrix} \\[2mm]
\sigma_1^{(0)}(\ba) &= \tau_1(\ba) :=  
\begin{pmatrix}
a_0+1-b_1, & a_1+1-b_1, & a_2+1-b_1 \\ 
           & 2-b_1,     & b_2+1-b_1
\end{pmatrix} \\
\sigma_2 ^{(0)}(\ba) &= \tau_2(\ba) :=  
\begin{pmatrix}
a_0+1-b_2, & a_1+1-b_2, & a_2+1-b_2 \\ 
           & b_1+1-b_2, & 2-b_2 
\end{pmatrix} \\[2mm]
\sigma_0 ^{(\infty)}(\ba) &= \sigma_0(\ba) :=  
\begin{pmatrix} 
a_0, & a_0+1-b_1, & a_0+1-b_2 \\
     & a_0+1-a_1, & a_0+1-a_2
\end{pmatrix} \\[2mm]
\sigma_1 ^{(\infty)}(\ba) &= \sigma_1(\ba) :=  
\begin{pmatrix}
a_1+1-b_1, & a_1, & a_1+1-b_2 \\                                    
     & a_1+1-a_0, & a_1+1-a_2
\end{pmatrix} \\[2mm] 
\sigma_2 ^{(\infty)}(\ba) &= \sigma_2(\ba) :=  
\begin{pmatrix}
a_2+1-b_2, & a_2+1-b_1, & a_2 \\                   
     & a_2+1-a_1, & a_2+1-a_0
\end{pmatrix}
\end{align*}
\caption{Six parameter involutions (including identity).} 
\label{tab:inv}
\end{table}
To describe other solutions to equation \eqref{eqn:ode} we introduce 
six parameter involutions $\sigma_i^{(\nu)}$, $i = 0,1,2$, 
$\nu = 0, \infty$, as in Table \ref{tab:inv} and impose condition   
\eqref{eqn:wd} for all of them, which becomes   
\begin{equation} \label{eqn:para-con}
a_i, \,\, a_i - a_j \not \in \Z, \qquad 
i, j = 0,1,2,3,4, \,\, i \neq j. 
\end{equation}
Let $\cS^{(\nu)}(\ba)$ be the space of local solutions to equation 
\eqref{eqn:ode} around $z = \nu \in \{0, \infty\}$.    
Under condition \eqref{eqn:para-con} none of the local exponent 
differences at $z = 0$ is an integer, so the local monodromy 
operator $\rM^{(0)} : \cS^{(0)}(\ba) \carl$ has three distinct 
eigenvalues $\lambda_i^{(0)}(\ba) :=\exp(-2 \ri \, \pi \, b_i)$, 
$i = 0,1, 2$, where $\ri := \sqrt{-1}$, and the corresponding 
eigen-solutions to \eqref{eqn:ode} are given by   
\begin{equation} \label{eqn:sol-zero}
y_i^{(0)}(\ba; z) := z^{1-b_i} \, \hgf(\sigma_i^{(0)}(\ba); z),  
\qquad i = 0, 1, 2. 
\end{equation}
Similarly, under condition \eqref{eqn:para-con} none of the local 
exponent differences at $z = \infty$ is an integer, so the local 
monodromy operator $\rM^{(\infty)} : \cS^{(\infty)}(\ba) \carl$ 
has three distinct eigenvalues $\lambda_i^{(\infty)}(\ba) := 
\exp(2 \ri \, \pi \, a_i)$, $i = 0,1, 2$, and the corresponding 
eigen-solutions to \eqref{eqn:ode} are given by    
\begin{equation} \label{eqn:sol-inf}
y_i^{(\infty)}(\ba; z) := e^{\mathrm{i} \pi  s(\sba) } 
\, z^{-a_i} \, \hgf(\sigma_i^{(\infty)}(\ba); 1/z),  
\qquad i = 0, 1, 2,   
\end{equation}
where $s(\ba)$ is defined in \eqref{eqn:rs}; the constant factor 
$\exp\left( \mathrm{i} \, \pi \, s(\ba) \right)$ plays no role in differential 
equation \eqref{eqn:ode}, but it will be important later when we discuss 
contiguous operators (see Lemma \ref{lem:contig-op}). 
As a summary the six functions in 
\eqref{eqn:sol-zero}-\eqref{eqn:sol-inf} are characterized  
by the eigen-relations 
\begin{equation} \label{eqn:eigen-r}
\rM^{(\nu)} \, y_i^{(\nu)}(\ba; z) = \lambda_i^{(\nu)}(\ba) 
\,  y_i^{(\nu)}(\ba; z) \qquad i = 0,1,2, \,\, \nu = 0, \infty, 
\end{equation}
uniquely up to nonzero constant multiples, where 
$\lambda_i^{(\nu)}(\ba)$ is invariant under any $\Z^5$-shift of 
$\ba$. 
\par  
Our main concern in this article is the hypergeometric series 
$\hgf(\ba) := \hgf(\ba; 1)$ of unit argument $z = 1$ as well as 
its six companions (including itself):  
\begin{equation} \label{eqn:companion}
y_i^{(\nu)}(\ba) := y_i^{(\nu)}(\ba; 1), \qquad i = 0, 1, 2, \quad 
\nu = 0, \infty,  
\end{equation}
where $\hgf(\ba) = y_0^{(0)}(\ba)$.  
It is well known that $\hgf(\ba)$ is convergent if and only if 
\begin{equation} \label{eqn:conv}
\rRe \, s(\ba) > 0,  
\end{equation} 
in which case the convergence is absolute and uniform in any compact 
subset so that the series $\hgf(\ba)$ defines a holomorphic function in 
domain \eqref{eqn:conv} with genericness condition \eqref{eqn:wd}.   
Observe that the Saalsch\"{u}tzian index is invariant under the 
six parameter changes in Table \ref{tab:inv},  
\begin{equation*} 
s(\ba) = s(\sigma_i^{(\nu)}(\ba)),  
\qquad i = 0,1,2, \,\, \nu = 0, \infty,  
\end{equation*} 
so all the series in \eqref{eqn:companion} are simultaneously convergent 
under the single condition \eqref{eqn:conv} together with genericness   
condition \eqref{eqn:para-con}.  
We shall denote by $\hgh(\ba; z)$ any function in 
\eqref{eqn:sol-zero}-\eqref{eqn:sol-inf} and by $\hgh(\ba)$ 
any function in \eqref{eqn:companion} respectively. 
The aim of this article is to establish the following.  
\begin{theorem} \label{thm:3tr}
Let $\hgh(\ba)$ be any member of the six functions in \eqref{eqn:companion}. 
If $\bp$, $\bq \in \Z^5$ are distinct shift vectors then the 
following assertions hold true:   
\begin{enumerate}
\item $\hgh(\ba+\bp)$ and $\hgh(\ba+\bq)$ are linearly independent 
over the field $\C(\ba)$.     
\item There exist unique rational functions $u(\ba)$, 
$v(\ba) \in \Q(\ba)$ such that 
\begin{equation} \label{eqn:3tr}
\hgh(\ba) = u(\ba) \, \hgh(\ba+\bp) + v(\ba) \, \hgh(\ba+\bq).   
\end{equation}
\item The rational functions $u(\ba)$ and $v(\ba)$ are common for  
all choices of $\hgh(\ba)$, that is, equation \eqref{eqn:3tr} is a 
simultaneous three-term relation for all functions in \eqref{eqn:companion}.    
\item There is a systematic recipe to determine $u(\ba)$ and 
$v(\ba)$ in finite steps $(\mathrm{Recipe} \, \ref{recipe})$. 
\item Three-term relation \eqref{eqn:3tr} admits 
an $S_2 \ltimes (S_3 \times S_3)$-symmetry of order seventy-two.    
\end{enumerate}
\end{theorem}
\par
Three-term relation \eqref{eqn:3tr} is interesting only when $\bp$ 
and $\bq$ are both different from zero, but the theorem itself is 
of course true even when one of them is zero. 
Obviously, for any triple of mutually distinct vectors $\bp$, $\bq$, 
$\br \in \Z^5$, we can  construct a general three-term relation 
\begin{equation} \label{eqn:g3tr} 
u(\ba) \, \hgh(\ba+\bp) + v(\ba) \, \hgh(\ba+\bq) + w(\ba) \, 
\hgh(\ba+\br) = 0,  
\end{equation}
having polynomial coefficients $u(\ba)$, $v(\ba)$, $w(\ba) \in \Q[\ba]$ 
without common factors, from a three-term relation of the form 
\eqref{eqn:3tr} by multiplying it by the common denominator of its 
coefficients and translating $\ba$ by an integer vector; 
\eqref{eqn:g3tr} is unique up to a constant multiple.  
\begin{remark} \label{rem:3tr}
If $K \supset \Q$ is a field of functions in $\ba = (a_0,a_1,a_2;a_3,a_4)$ 
such that 
\begin{equation} \label{eqn:transc}
\mbox{$a_0$, $a_1$, $a_2$, $a_3$, $a_4$ are algebraically independent 
over $K$}, 
\end{equation}
then assertion (1) of Theorem \ref{thm:3tr} remains true over the field 
$K(\ba)$. 
For example this condition is fulfilled by the field $K = \cM_{\vL}$ of 
$\vL$-periodic meromorphic functions on $\C^5$ for a lattice $\vL$, 
where a {\sl lattice} in $\C^5$ is the $\Z$-linear span of five 
vectors that form a $\C$-linear basis of $\C^5$.     
\end{remark}
\par
This remark is of particular interest because the six functions in 
\eqref{eqn:companion} satisfy some linear relations over 
the periodic function field $\cM_{\vL}$ with period 
lattice $\vL = (2 \Z)^5$ but not over the field $\C$.    
For example the three functions associated with $z = 0$ satisfy a 
linear relation    
\begin{equation} \label{eqn:mimachi}
c_0(\ba) \, y_0^{(0)}(\ba) + c_1(\ba) \, y_1^{(0)}(\ba) +c_2(\ba) \, y_2^{(0)}(\ba) = 0,  
\end{equation}
where the coefficients $c_i(\ba) \in \cM_{\vL}$ are given in terms of 
trigonometric functions by 
\[
c_i(\ba) := c(\sigma_i^{(0)}(\ba)) \qquad \mbox{with} \quad  
c(\ba) := \dfrac{\sin \pi a_0 \cdot \sin \pi a_1 \cdot 
\sin \pi a_2}{\sin \pi b_1 \cdot \sin \pi b_2}.   
\]
One way to prove relation \eqref{eqn:mimachi} is to take the difference of  
two connection formulas  
\begin{equation} \label{eqn:Thomae}
e^{-\mathrm{i} \pi s(\sba)} y_0^{(\infty)}(\ba) = 
\dfrac{ \sin \pi a_1 \cdot \sin \pi a_2}{\sin \pi b_j  
\cdot \sin \pi (b_k -a_0)} \, y_0^{(0)}(\ba) 
- 
\dfrac{\sin \pi(a_1-b_j) \cdot \sin \pi(a_2-b_j)}{\sin \pi b_j 
\cdot \sin \pi(b_k - a_0)} \, y_j^{(0)}(\ba),  
\end{equation} 
with $(j, k) = (1, 2)$, $(2, 1)$.   
Essentially, either formula in \eqref{eqn:Thomae} is Thomae's second 
fundamental relation in disguise; thanks to the normalization \eqref{eqn:hgf(z)} 
and definition \eqref{eqn:companion}, it gains a better appearance than 
its original form in Thomae \cite[page 72]{Thomae} (see also 
Bailey \cite[\S 3.2, formula (2)]{Bailey1}) for which the 
periodic nature of its coefficients is not transparent.  
\par
In Theorem \ref{thm:3tr} the symmetry on three-term relations 
in assertion (5) is just a corollary to the uniqueness and 
simultaneousness in assertions (2) and (3), but this issue 
is interesting in its own light and hence discussed in 
\S \ref{sec:symmetry}, which contains a detailed account  
of assertion (5).       
\section{Differential Operators} \label{sec:diff-op}
We take the differential operator approach to the contiguous relations 
as in Bailey \cite{Bailey2}, but we are more systematic with this 
method and keep it minimum for the sake of conciseness. 
Let 
\begin{alignat*}{3}
\be_0 &:= (1,0,0;0,0), \qquad & \be_1 &:= (0,1,0;0,0), \qquad & \be_2 &:= (0,0,1;0,0), \\ 
\be_3 &:= (0,0,0;1,0), \qquad & \be_4 &:= (0,0,0;0,1), \qquad & \1    &:= (1,1,1;1,1).       
\end{alignat*}
Our approach to the $\hgF$ contiguous relations is based on the following.  
\begin{lemma} \label{lem:contig-op} 
For any member $\hgh(\ba; z)$ of the six functions in 
\eqref{eqn:sol-zero}-\eqref{eqn:sol-inf}, we have 
\begin{subequations} \label{eqn:contig-op}
\begin{alignat}{2}
(\theta + a_i) \, \hgh(\ba; z) &= \hgh(\ba + \be_i; z) \qquad & 
i &= 0, 1, 2, \label{eqn:contig-op1} \\
(\theta + a_i -1) \, \hgh(\ba; z) &= \hgh(\ba - \be_i; z) \qquad & 
i &= 3, 4, \label{eqn:contig-op2} \\
\partial \, \hgh(\ba; z) &= \hgh(\ba + \1; z) &  
\partial &:= \ts \frac{\partial}{\partial z}. \label{eqn:contig-op3}
\end{alignat}
\end{subequations}
\end{lemma}
{\it Proof}.  
Using the basic commutation relations $\theta \cdot z = z \cdot 
(\theta+1)$ and $\theta \cdot \partial = \partial \cdot (\theta-1)$, 
we can easily observe that the differential operator $L(\ba)$ in 
\eqref{eqn:ode} admits commutation relations 
\begin{alignat*}{2}
L(\ba + \be_i) \, (\theta + a_i) &= (\theta + a_i) \, L(\ba)  
\qquad & i &= 0, 1, 2,  \\
L(\ba - \be_i) \, (\theta + a_i-1) &= (\theta + a_i-2) \, L(\ba)  
\qquad & i &= 3, 4, \\
z \, L(\ba + \1) \, \partial &= (\theta -1) \, L(\ba). & & 
\end{alignat*}
This implies that for any $j = 0,1,2$, $\nu = 0, \infty$ 
we have on the one hand, 
\begin{alignat*}{2}
(\theta + a_i) \, y_j^{(\nu)}(\ba; z) &\in \cS^{(\nu)}(\ba + \be_i; z) \qquad & i &= 0, 1, 2, \\
(\theta + a_i-1) \, y_j^{(\nu)}(\ba; z) &\in \cS^{(\nu)}(\ba - \be_i; z) \qquad & i &= 3, 4, \\
\partial \, y_j^{(\nu)}(\ba; z) &\in \cS^{(\nu)}(\ba +\1; z). & &   
\end{alignat*}
On the other hand, since the local monodromy operator $\rM^{(\nu)}$ is commutative with 
differential operators with rational coefficients in $z$ and the eigenvalues 
$\lambda_j^{(\nu)}(\ba)$ are $\Z^5$-periodic in $\ba$,   
\begin{alignat*}{2}
\rM^{(\nu)} \cdot (\theta + a_i) \, y_j^{(\nu)}(\ba; z) &= 
\lambda_j^{(\nu)}(\ba+ \be_i) \cdot (\theta + a_i) \, y_j^{(\nu)}(\ba; z) \qquad & i &= 0, 1, 2, \\
\rM^{(\nu)} \cdot (\theta + a_i-1) \, y_j^{(\nu)}(\ba; z) &= 
\lambda_j^{(\nu)}(\ba- \be_i) \cdot (\theta + a_i-1)  \, y_j^{(\nu)}(\ba; z) \qquad & i &= 3, 4, \\
\rM^{(\nu)} \cdot \partial \, y_j^{(\nu)}(\ba; z) &= 
\lambda_j^{(\nu)}(\ba + \1) \cdot \partial \, y_j^{(\nu)}(\ba; z). & &   
\end{alignat*}
Since eigen-solutions \eqref{eqn:sol-zero}-\eqref{eqn:sol-inf} are uniquely 
determined by eigen-relations \eqref{eqn:eigen-r} up to constant multiples, 
there must be constants $\alpha_{ij}^{(\nu)}(\ba)$ and 
$\beta_j^{(\nu)}(\ba)$ with respect to $z$ such that 
\begin{alignat}{2}
(\theta + a_i) \, y_j^{(\nu)}(\ba; z) &= 
\alpha_{ij}^{(\nu)}(\ba) \cdot y_j^{(\nu)}(\ba + \be_i; z) \qquad & 
i &= 0, 1, 2,  \nonumber \\
(\theta + a_i-1) \, y_j^{(\nu)}(\ba; z) &=  
\alpha_{ij}^{(\nu)}(\ba) \cdot y_j^{(\nu)}(\ba - \be_i; z) \qquad & 
i &= 3, 4, \label{eqn:const} \\
\partial \, y_j^{(\nu)}(\ba; z) &= 
\beta_j^{(\nu)}(\ba) \cdot y_j^{(\nu)}(\ba +\1; z). & &  \nonumber  
\end{alignat}
These constants can be determined by equating the leading coefficients of both sides 
in \eqref{eqn:const}. 
Thanks to the normalization in \eqref{eqn:hgf(z)} and the constant factor 
$\exp(\ri \, \pi \, s(\ba))$ in \eqref{eqn:sol-inf} we have 
\[
\alpha_{ij}^{(\nu)}(\ba) = \beta_{j}^{(\nu)}(\ba) = 1 \qquad 
i = 0,1,2,3,4, \,\, j = 0,1,2, \,\, \nu = 0, \infty  
\quad (\mbox{simultaneously one !}), 
\]
which together with \eqref{eqn:const} yields the desired formulas  
\eqref{eqn:contig-op}. \hfill $\Box$ \par\medskip 
The differential equation \eqref{eqn:ode} is written in terms of $\theta$. 
Expressed in $\partial$, it is equivalent to   
\begin{equation} \label{eqn:ode2}
\begin{split}
M(\ba) \, y &:= \big[ (1-z) \, z^2 \, \partial^3 + 
z \, \{ \, ( a_3 + a_4 + 1) -( \varphi_1(\ba) + 3 ) \, z \, \} \, \partial^2 \\
&\qquad + \{\, a_3 a_4 - ( \, \varphi_2(\ba) + \varphi_1(\ba) +1 \,) \, z \, \} 
\, \partial - \varphi_3(\ba) \big] \, y = 0,  
\end{split}
\end{equation}
where $\varphi_i(\ba)$ is the $i$-th elementary symmetric polynomial in 
the numerator parameters $(a_0, a_1, a_2)$:  
\[
\varphi_1(\ba) := a_0 + a_1 + a_2, \qquad  
\varphi_2(\ba) := a_0 a_1 + a_1 a_2 + a_2 a_0, \qquad   
\varphi_3(\ba) := a_0 a_1 a_2.  
\]
The differential equation \eqref{eqn:ode2} leads to an important 
three-term relation 
\begin{equation} \label{eqn:one-two}
\{s(\ba) -2 \} \, \hgh(\ba +\2) + \psi(\ba) \, \hgh(\ba +\1) - 
\varphi_3(\ba) \, \hgh(\ba) = 0,    
\end{equation}
with $\2 := (2,2,2;2,2)$ and $\psi(\ba) := 
a_3 a_4 -\varphi_2(\ba) - \varphi_1(\ba) -1$. 
Indeed, substitute $y = \hgh(\ba; z)$ into equation \eqref{eqn:ode2}, 
use formula \eqref{eqn:contig-op3} and put $z = 1$ to get 
equation \eqref{eqn:one-two}. 
It indicates explicitly how $\hgh(\ba+\2)$ can be written as a linear 
combination of $\hgh(\ba)$ and $\hgh(\ba+\1)$. 
\par
For $i = 0,1,2$ with $\{i,j,k\} = \{0,1,2\}$, one has    
\begin{subequations} \label{eqn:crnp}
\begin{align}
\hgh(\ba+\be_i)    &= a_i \, \hgh(\ba) + \hgh(\ba+\1),  
\label{eqn:crnp-a} \\
\hgh(\ba+\be_i+\1) &= \dfrac{\varphi_3(\ba)}{s(\ba) -2} \, \hgh(\ba) + 
\dfrac{a_j a_k -(a_i-a_3+1)(a_i-a_4+1)}{s(\ba) -2} \, \hgh(\ba+\1), 
\label{eqn:crnp-b}
\end{align} 
\end{subequations}
while for $i = 3,4$, one has   
\begin{subequations} \label{eqn:crdm}
\begin{align}
\hgh(\ba-\be_i) &= (a_i -1) \, \hgh(\ba) + \hgh(\ba+\1),  
\label{eqn:crdm-a} \\
\hgh(\ba-\be_i+\1) &= \dfrac{\varphi_3(\ba)}{s(\ba) -2} \, \hgh(\ba) +  
\dfrac{(a_i-1)^2-\varphi_1(\ba) \, (a_i-1) + \varphi_2(\ba)}{s(\ba) -2} \, 
\hgh(\ba+\1). 
\label{eqn:crdm-b}
\end{align} 
\end{subequations}  
Indeed, formula \eqref{eqn:crnp-a} resp. \eqref{eqn:crdm-a} is  
obtained by using \eqref{eqn:contig-op3} in \eqref{eqn:contig-op1} 
resp. \eqref{eqn:contig-op2} and putting $z = 1$. 
Formulas \eqref{eqn:crnp-b} resp. \eqref{eqn:crdm-b} is then   
obtained by replacing $\ba$ with $\ba+\1$ in \eqref{eqn:crnp-a} resp.   
\eqref{eqn:crdm-a} and eliminating $\hgh(\ba+\2)$ through three-term 
relation \eqref{eqn:one-two}.  
\section{Contiguous and Connection Matrices} \label{sec:cm}
To establish our main theorem it is crucial to put contiguous relations 
into matrix forms. 
If 
\begin{equation} \label{eqn:vec-h}
\bhgh(\ba) := 
\begin{pmatrix} \hgh(\ba) \\ \hgh(\ba + \1) \end{pmatrix},  
\end{equation}
then the matrix version of contiguous relations is represented as 
\begin{equation} \label{eqn:contig-m}
\bhgh(\ba + \varepsilon \, \be_i) = A_i^{\varepsilon}(\ba) \, \bhgh(\ba)
\qquad i = 0,1,2,3,4, \,\, \ve = \pm,   
\end{equation}
where the matrices $A_i^{\varepsilon}(\ba)$ are given in 
\eqref{eqn:cm} of Table \ref{tab:cm}; they are invertible over the 
field $\Q(\ba)$, having non-vanishing determinants as in 
\eqref{eqn:det-cm}.    
\begin{table}[pp]
\begin{subequations} \label{eqn:cm}
\begin{alignat}{2}
A_i^+(\ba) &= 
\begin{pmatrix}
a_i & 1 \\[2mm]
\dfrac{\varphi_3(\ba)}{s(\ba) -2} & 
\dfrac{a_j a_k -(a_i-a_3+1)(a_i-a_4+1)}{s(\ba) -2}
\end{pmatrix} 
\quad & i &= 0, 1, 2, \label{eqn:cm1} \\[3mm]  
A_i^+(\ba) &= 
\begin{pmatrix}
\dfrac{a_i^2 - \varphi_1(\ba) \, a_i + \varphi_2(\ba)}{(a_i-a_0)(a_i-a_1)(a_i-a_2)} & 
- \dfrac{s(\ba)-1}{(a_i-a_0)(a_i-a_1)(a_i-a_2)} \\[4mm]
- \dfrac{\varphi_3(\ba)}{(a_i-a_0)(a_i-a_1)(a_i-a_2)} & 
\dfrac{a_i \, \{s(\ba)-1\}}{(a_i-a_0)(a_i-a_1)(a_i-a_2)} 
\end{pmatrix} 
\quad & i &= 3, 4, \label{eqn:cm2} \\[3mm]
A_i^-(\ba) &= 
\begin{pmatrix}
\dfrac{(a_i-a_3)(a_i-a_4) - a_j a_k}{(a_i-1)(a_i-a_3)(a_i-a_4)} & 
\dfrac{s(\ba)-1}{(a_i-1)(a_i-a_3)(a_i-a_4)} \\[4mm]
\dfrac{a_j a_k}{(a_i-a_3)(a_i-a_4)} & - \dfrac{s(\ba)-1}{(a_i-a_3)(a_i-a_4)} 
\end{pmatrix} \quad & i &= 0, 1, 2, \label{eqn:cm3}  \\[3mm]
A_i^-(\ba) &= 
\begin{pmatrix}
a_i -1 & 1 \\[2mm]
\dfrac{\varphi_3(\ba)}{s(\ba) -2} & 
\dfrac{(a_i-1)^2-\varphi_1(\ba) \, (a_i-1) + \varphi_2(\ba)}{s(\ba) -2}
\end{pmatrix} 
\quad & i &= 3, 4. \label{eqn:cm4} 
\end{alignat} 
\end{subequations}
\begin{subequations} \label{eqn:det-cm}
\begin{alignat}{2}
\det A_i^+(\ba) &= - \dfrac{a_i(a_i-a_3+1)(a_i-a_4+1)}{s(\ba) -2} \qquad & 
i &= 0, 1, 2,  \label{eqn:det-cm1}  \\[2mm]
\det A_i^+(\ba) &= \dfrac{s(\ba)-1}{(a_i - a_0)(a_i - a_1)(a_i - a_2)} \qquad & 
i &= 3, 4,  \label{eqn:det-cm2}  \\[2mm]
\det A_i^-(\ba) &= - \dfrac{s(\ba) -1}{(a_i - 1)(a_i - a_3)(a_i - a_4)} \qquad & 
i &= 0, 1, 2,  \label{eqn:det-cm3} \\[2mm]
\det A_i^-(\ba) &= \dfrac{(a_i-a_0-1)(a_i-a_1-1)(a_i-a_2-1)}{s(\ba) -2} \qquad & 
i &= 3,4,  \label{eqn:det-cm4} 
\end{alignat}
\end{subequations}
\caption{Contiguous matrices and their determinants; 
$\{i,j,k\} = \{0,1,2\}$ for $i = 0, 1, 2$.} \label{tab:cm}
\begin{subequations} \label{eqn:pp}
\begin{alignat}{2}
B_i(\ba) &= 
\begin{pmatrix}
a_i & 1 \\[2mm]
\dfrac{\varphi_3(\ba)}{s(\ba)} & 
\dfrac{a_j a_k -(a_i-a_3)(a_i-a_4)}{s(\ba)}
\end{pmatrix} 
\quad & i &= 0, 1, 2, \label{eqn:pp1} \\[3mm]  
B_i(\ba) &= 
\begin{pmatrix}
\dfrac{a_i^2 - \varphi_1(\ba) \, a_i + \varphi_2(\ba)}{(a_i-a_0)(a_i-a_1)(a_i-a_2)} & 
- \dfrac{s(\ba)}{(a_i-a_0)(a_i-a_1)(a_i-a_2)} \\[4mm]
- \dfrac{\varphi_3(\ba)}{(a_i-a_0)(a_i-a_1)(a_i-a_2)} & 
\dfrac{a_i \, s(\ba) }{(a_i-a_0)(a_i-a_1)(a_i-a_2)} 
\end{pmatrix} 
\quad & i &= 3, 4, \label{eqn:pp2} 
\end{alignat} 
\end{subequations}
\caption{Principal parts of contiguous matrices; 
$\{i,j,k\} = \{0,1,2\}$ for $i = 0, 1, 2$.} 
\label{tab:pp}
\end{table}
Formulas \eqref{eqn:cm1} and \eqref{eqn:cm4} come directly from 
\eqref{eqn:crnp} and \eqref{eqn:crdm}, while the remaining formulas 
\eqref{eqn:cm2} and \eqref{eqn:cm3} are derived from them via  
the relation $A_i^{\varepsilon}(\ba) = 
A_i^{-\varepsilon}(\ba + \varepsilon \, \be_i)^{-1}$. 
The contiguous matrices satisfy compatibility conditions:  
\begin{equation} \label{eqn:cpb}
A_i^{\ve_i}(\ba + \ve_j \, \be_j) \, A_j^{\ve_j}(\ba) = 
A_j^{\ve_j}(\ba + \ve_i \, \be_i) \, A_i^{\ve_i}(\ba), 
\qquad i,j = 0,1,2,3,4, \quad \ve_i, \ve_j = \pm. 
\end{equation} 
\par
Given an integer vector $\bp \in \Z^5$, a finite sequence 
$\bp_0, \, \bp_1, \, \dots, \, \bp_m$ of integer vectors in $\Z^5$ is said 
to be a {\sl lattice path} from the origin $\0$ to $\bp$ if there exist an 
$m$-tuple of indices $(i_1, \dots, i_m) \in \{0,1,2,3,4\}^m$ and an 
$m$-tuple of signs $(\ve_1, \dots, \ve_m) \in \{ \pm \}^m$ such that 
\begin{equation} \label{eqn:pm}
\bp_0 = \0, \qquad \bp_m = \bp, \qquad 
\bp_k - \bp_{k-1} = \ve_k \, \be_{i_k}, \qquad k = 1, \dots, m.  
\end{equation}
To such a lattice path one can associate the products of contiguous 
matrices 
\begin{equation} \label{eqn:A(a;p)} 
A(\ba; \bp) := A_{i_m}^{\ve_m}(\ba+\bp_{m-1}) \, 
A_{i_{m-1}}^{\ve_{m-1}}(\ba+\bp_{m-2}) \, \cdots \, 
A_{i_2}^{\ve_2}(\ba+\bp_1) \, A_{i_1}^{\ve_1}(\ba+\bp_0).  
\end{equation}
Thanks to the compatibility condition \eqref{eqn:cpb} the ensuing matrix 
$A(\ba; \bp)$ depends only on the terminal vector $\bp$, being independent 
of which lattice path is chosen. 
In view of \eqref{eqn:contig-m} and \eqref{eqn:A(a;p)} the effect 
of translation $\ba \mapsto \ba + \bp$ on the function $\bhgh(\ba)$ 
is represented by   
\begin{equation} \label{eqn:shift-m}
\bhgh(\ba+\bp) = A(\ba; \bp) \, \bhgh(\ba).   
\end{equation}
We refer to $A(\ba;\bp)$ as the {\sl connection matrix} for a given  
shift vector $\bp \in \Z^5$. 
Induction on the length $m$ in \eqref{eqn:A(a;p)} based on the determinant 
formulas \eqref{eqn:det-cm} yields 
\begin{equation} \label{eqn:det}
\det A(\ba; \bp) = \delta(\ba; \bp) :=  
\dfrac{ (-1)^{p_0+p_1+p_2} \cdot (s(\ba)-1, \, s(\bp) ) \cdot 
\prod_{i=0}^2 (a_i, \, p_i)}{\prod_{i=0}^2 \prod_{j=3}^4 (a_j-a_i, \, p_j-p_i)} \neq 0 
\quad \mbox{in} \quad \Q(\ba),   
\end{equation}
where $\bp = (p_0, p_1, p_2; p_3, p_4)$ and  $s(\bp)$ is its Saalsch\"{u}tzian index. 
\par
Focusing on the upper row of connection matrix $A(\ba; \bp)$, if we write 
\begin{equation} \label{eqn:ur}
A(\ba; \bp) = \begin{pmatrix} r_1(\ba; \bp)  & r(\ba; \bp) \\ * & * 
\end{pmatrix},   
\end{equation} 
then the upper entry of matrix equation \eqref{eqn:shift-m} leads to a three-term 
relation 
\begin{equation} \label{eqn:one-p}
\hgh(\ba + \bp) = r_1(\ba; \bp) \, \hgh(\ba) + r(\ba; \bp) \, \hgh(\ba + \1). 
\end{equation}
Similarly, for another vector $\bq \in \Z^5$ different from $\bp$ one has   
a second three-term relation  
\begin{equation} \label{eqn:one-q}
\hgh(\ba + \bq) = r_1(\ba; \bq) \, \hgh(\ba) + r(\ba; \bq) \, \hgh(\ba + \1).  
\end{equation}
\par
The linear system \eqref{eqn:one-p}-\eqref{eqn:one-q} can be written in the 
matrix form  
\begin{equation} \label{eqn:base-change}
\begin{pmatrix}
\hgh(\ba+\bp) \\
\hgh(\ba+\bq) 
\end{pmatrix}
= 
\begin{pmatrix}
r_1(\ba;\bp) & r(\ba;\bp) \\
r_1(\ba;\bq) & r(\ba;\bq) 
\end{pmatrix}
\begin{pmatrix}
\hgh(\ba) \\ 
\hgh(\ba+\1) 
\end{pmatrix}, 
\end{equation}
and we wonder whether the system is invertible, that is, if the determinant 
\begin{equation} \label{eqn:nonzero}
\vD(\ba;\bp,\bq) :=
\left|\,  
\begin{matrix}
r_1(\ba; \bp) & r(\ba; \bp) \\ 
r_1(\ba; \bq) & r(\ba; \bq) 
\end{matrix}\, \right|  
\end{equation}
is non-vanishing in $\Q(\ba)$. 
To answer this question we first prove the following lemma. 
\begin{lemma} \label{lem:solve} 
With formula \eqref{eqn:det} and definition \eqref{eqn:nonzero} one has 
\begin{equation} \label{eqn:solve}
\vD(\ba; \bp, \bq) =  \det A(\ba; \bp) \cdot r(\ba+\bp; \bq-\bp)
= \delta(\ba; \bp) \cdot r(\ba+\bp; \bq-\bp).  
\end{equation}
\end{lemma}
{\it Proof}. By the chain rule $A(\ba; \bq) = A(\ba+\bp; \bq-\bp) A(\ba; \bp)$ 
one has   
\begin{align*}
A(\ba+\bp; \bq-\bp) &= A(\ba; \bq) A(\ba; \bp)^{-1} 
= \begin{pmatrix} r_1(\ba; \bq) & r(\ba; \bq) \\ * & * \end{pmatrix} 
\dfrac{1}{\det A(\ba; \bp) } 
 \begin{pmatrix} * & -r(\ba; \bp) \\ * &  r_1(\ba; \bp) \end{pmatrix} \\[2mm]
&= \dfrac{1}{\det A(\ba; \bp) } 
\begin{pmatrix} * & \vD(\ba; \bp, \bq) \\ * &  * \end{pmatrix},  
\end{align*}
the $(1, 2)$ entry of which leads to $r(\ba+\bp; \bq-\bp) = 
\vD(\ba; \bp, \bq)/ \det A(\ba; \bp)$. \hfill $\Box$ \par \medskip
Thus the vanishing of $\vD(\ba; \bp, \bq)$ in $\Q(\ba)$ is equivalent to the 
condition that the matrix $A(\ba+\bp; \bq-\bp)$ or more simply 
$A(\ba; \bq-\bp)$ be lower triangular over $\Q(\ba)$.  
In the next section we show that this never happens unless 
$\bp = \bq$ (see Proposition \ref{prop:lt}).      
\section{Principal Parts and Triangular Matrices} \label{sec:pp-tm}
The non-commutative nature of contiguous matrices makes the structure of 
$A(\ba; \bp)$ so intricate that a direct tackle to our problem may be difficult.     
In such a situation, taking the ``principal parts" of non-commutative 
objects often turns things commutative, although much information is 
usually lost in this process, but the ensuing commutative objects, 
perhaps much easier to handle, may still retain something important. 
Such an idea is employed in Iwasaki \cite[\S 9]{Iwasaki} 
to analyze contiguous matrices of Gauss's hypergeometric functions  
${}_2F_1$. 
\par
We now apply a similar thought to $\hgF(1)$, but with a necessary twist 
involving 
\[
D(t) := \mathrm{diag}\{1, \, t \}, \qquad 
\delta_i := 
\begin{cases} 
-1 & (i = 0,1,2), \\ 
+1 & (i = 3,4),  
\end{cases}
\]
where $t$ is a formal parameter. 
An inspection shows that there is a formal Laurent expansion    
\begin{equation} \label{eqn:pp-cm}
t^{\delta_i \ve} \cdot D(t)^{-1} \, A_i^{\varepsilon}(t \, \ba) \, D(t) 
= B_i(\ba)^{\ve} + O(1/t) \qquad \mbox{around} \quad t = \infty, 
\end{equation}
where $B_i(\ba)^{\ve} := B_i(\ba)^{\pm 1}$ for $\ve = \pm$ and 
$B_i(\ba)$ is given as in \eqref{eqn:pp} of Table \ref{tab:pp}. 
The matrix $B_i(\ba)$ is obtained from $A_i^+(\ba)$ by taking the top 
homogeneous components of its entries, where the top homogeneous 
component of a rational function is the ratio of the top homogeneous 
polynomials of its numerator and denominator.   
Observe that     
\begin{subequations} \label{eqn:det-pp}
\begin{alignat}{2}
\det B_i(\ba) &= - \dfrac{a_i(a_i-a_3)(a_i-a_4)}{s(\ba)} \qquad & 
i &= 0, 1, 2,  \label{eqn:det-pp1}  \\[2mm]
\det B_i(\ba) &= \dfrac{s(\ba)}{(a_i - a_0)(a_i - a_1)(a_i - a_2)} \qquad & 
i &= 3, 4,  \label{eqn:det-pp2} 
\end{alignat}
\end{subequations} 
so the matrices $B_i(\ba)$ are invertible over $\Q(\ba)$. 
Compatibility condition \eqref{eqn:cpb} for contiguous matrices 
naturally leads to the commutativity of their principal parts: 
\begin{equation} \label{eqn:comm} 
B_i(\ba) B_j(\ba) = B_j(\ba) B_i(\ba) \qquad i,j = 0,1,2,3,4. 
\end{equation} 
\begin{lemma} \label{lem:pp} 
For any $\bp = (p_0,p_1,p_2;p_3,p_4) \in \Z^5$ there is a formal Laurent 
expansion  
\begin{equation} \label{eqn:pp-A(a;p)}
t^{s(\sbp)} \cdot D(t)^{-1} \, A(t \, \ba; \bp) \, D(t) 
= B(\ba; \bp) + O(1/t) \qquad \mbox{around} \quad t = \infty, 
\end{equation}
where the matrix $B(\ba; \bp)$, called the principal part of 
$A(\ba; \bp)$, is given by 
\begin{equation} \label{eqn:B(a;p)}
B(\ba; \bp) = B_0(\ba)^{p_0} \, B_1(\ba)^{p_1} \, 
B_2(\ba)^{p_2} \, B_3(\ba)^{p_3} \, B_4(\ba)^{p_4}.  
\end{equation}
\end{lemma}
{\it Proof}. 
Replacing $\ba$ with $t \,\ba$ in formula \eqref{eqn:A(a;p)}, 
substituting the result into the left-hand side of 
\eqref{eqn:pp-A(a;p)} and using formula \eqref{eqn:pp-cm}, 
we observe that expansion \eqref{eqn:pp-A(a;p)} holds true with 
\[
B(\ba; \bp) = B_{i_m}(\ba)^{\ve_m} \, 
B_{i_{m-1}}(\ba)^{\ve_{m-1}} \, \cdots \, 
B_{i_2}(\ba)^{\ve_2} \, B_{i_1}(\ba)^{\ve_1}.  
\]
The commutativity \eqref{eqn:comm} allows us to rearrange the order 
of the matrix products so that 
\[
B(\ba; \bp) = 
\prod_{k \in \vL_0} B_{0}(\ba)^{\ve_k} \, 
\prod_{k \in \vL_1} B_{1}(\ba)^{\ve_k} \, 
\prod_{k \in \vL_2} B_{2}(\ba)^{\ve_k} \,
\prod_{k \in \vL_3} B_{3}(\ba)^{\ve_k} \,
\prod_{k \in \vL_4} B_{4}(\ba)^{\ve_k},    
\]
where $\vL_j := \{ k = 1, \dots, m \,:\, i_k = j \, \}$ for $j = 0,1,2,3,4$. 
Since $\ds \sum_{k \in \vL_j} \ve_k = p_j$ we have \eqref{eqn:B(a;p)}. 
\hfill $\Box$ 
\begin{remark} \label{rem:pp} 
Since $D(t)$ is a diagonal matrix, expansion formula \eqref{eqn:pp-A(a;p)} 
implies that if $A(\ba; \bp)$ is a lower triangular matrix over the field 
$\Q(\ba)$, then so must be $B(\ba; \bp)$. 
\end{remark}
\par 
The aim of this section is to prove the following proposition. 
\begin{proposition} \label{prop:lt} 
For any shift vector $\bp = (p_0,p_1,p_2;p_3,p_4) \in \Z^5$, if $A(\ba; \bp)$ is lower 
triangular over $\Q(\ba)$, that is, if $r(\ba; \bp)$ vanishes in $\Q(\ba)$ then 
$\bp$ must be $\0$.  
\end{proposition}
\par
By Remark \ref{rem:pp}, to prove the proposition it is sufficient to 
show the following lemma.  
\begin{lemma} \label{lem:lt} 
If $B(\ba; \bp)$ is lower triangular over $\Q(\ba)$ then $\bp$ 
must be $\0$.  
\end{lemma}
\par
Thanks to the commutative nature of principal part matrices, 
Lemma \ref{lem:lt} is much more tractable because simultaneous 
diagonalization technique is available. 
Our strategy to prove the lemma is {\sl specializations}, 
that is, to work with some special but convenient values of $\ba$. 
For this purpose we may take any complex value of $\ba$ 
as far as the matrices $B_i(\ba)$, $i = 0,1,2,3,4$, have no 
zero denominators or no zero determinants; such a value is 
called {\sl good}. 
In view of the denominators in Table \ref{tab:pp} and 
the determinants in \eqref{eqn:det-pp} the goodness condition  
is given by 
\begin{equation} \label{eqn:good} 
a_0 a_1 a_2 \cdot s(\ba) \cdot \prod_{i=1}^2 \prod_{j=0}^2 
(b_i - a_j) \neq 0. 
\end{equation} 
Here we do not have to mind the genericness condition \eqref{eqn:para-con}, 
which has nothing to do with our current concern. 
If Lemma \ref{lem:lt} is proved for some good values of $\ba$, 
then we are done. 
\[
\mbox{An invertible matrix 
$P = \begin{pmatrix} \alpha & \beta \\ \gamma & \delta \end{pmatrix}$ 
is said to be {\sl admissible} if $\alpha \, \beta \neq 0$}. 
\]
Our discussion is then based on the following simple lemma 
for lower triangular matrices.     
\begin{lemma} \label{lem:eev} 
If an admissible matrix $P$ diagonalizes a lower triangular 
matrix $B$ as $P^{-1} B P = \mathrm{diag}\{\lambda_1, \, \lambda_2\}$, 
then one must have $\lambda_1 = \lambda_2$ and so $B$ must be a scaler 
matrix.  
\end{lemma}
{\it Proof}. 
A straightforward calculation shows 
\[
B = P \begin{pmatrix} \lambda_1 & 0 \\ 0 & \lambda_2 \end{pmatrix} 
P^{-1} = \frac{1}{\varDelta} 
\begin{pmatrix} 
\lambda_1 \, \alpha \, \delta - \lambda_2 \, \beta \, \gamma & 
(\lambda_2 - \lambda_1) \, \alpha \, \beta \\ 
(\lambda_1-\lambda_2) \, \gamma \, \delta & 
\lambda_2 \, \alpha \, \delta - \lambda_1 \, \beta \, \gamma
\end{pmatrix}, \qquad \varDelta := \det P \neq 0.  
\]
Since $B$ is lower triangular, we have $(\lambda_2 - \lambda_1) \alpha \, \beta = 0$ 
and hence $\lambda_2 - \lambda_1 = 0$ by $\alpha \, \beta \neq 0$. 
\hfill $\Box$ \par\medskip
When working with a special value of $\ba$, we abbreviate 
$B(\ba; \bp)$ to $B$ and $B_i(\ba)$ to $B_i$.      
\begin{lemma} \label{lem:lt1} 
If $B(\ba; \bp)$ is lower triangular over $\Q(\ba)$ then
\begin{equation} \label{eqn:lt1} 
p_0+p_1+p_2 = p_3+p_4 = 0. 
\end{equation} 
\end{lemma}
{\it Proof}. 
Taking a specialization $\ba = (2,2,2;-7,-7)$, which is good from 
\eqref{eqn:good}, we have 
\[
B_0 = B_1 = B_2 = 
\begin{pmatrix} 2 & 1 \\ -2/5 & 77/ 20 \end{pmatrix}, \qquad 
B_3 = B_4 = 
\begin{pmatrix} 
-103/729 & -20/729 \\  8/729 & -140/729 
\end{pmatrix}. 
\]
These matrices are simultaneously diagonalized by an admissible matrix    
\[
P = 
\begin{pmatrix}
5/8 & 4 \\ 1 & 1 
\end{pmatrix}
\qquad \mbox{as} \qquad  
\begin{array}{rcl}
P^{-1} \, B_0 \, P &=& \diag \{ \, 2 \cdot 3^2 \cdot 5^{-1},  \quad 2^{-2} \cdot 3^2 \, \}, \\  
P^{-1} \, B_3 \, P &=& \diag \{ \, - 3^{-3} \cdot 5, \quad -2^2 \cdot 3^{-3} \, \}.  
\end{array}
\]
Formula \eqref{eqn:B(a;p)} then implies 
$B = B_0^{m} \, B_3^{n}$ with $m := p_0+p_1+p_2$ and $n := p_3+p_4$. 
So the lower triangular matrix $B$ is diagonalized as 
$P^{-1} \, B \, P = \diag\{\lambda_1, \, \lambda_2\}$, where  
\begin{align*}
\lambda_1 &= (2 \cdot 3^2 \cdot 5^{-1})^m ( - 3^{-3} \cdot 5)^n 
= (-1)^n \cdot 2^m \cdot 3^{2 m-3 n} \cdot 5^{n-m}, \\
\lambda_2 &= (2^{-2} \cdot 3^2)^m (- 2^2 \cdot 3^{-3} )^n 
= (-1)^n \cdot 2^{2 n - 2 m} \cdot 3^{2 m-3 n}.  
\end{align*} 
By Lemma \ref{lem:eev} one must have $\lambda_1 = \lambda_2$, that is, 
$\lambda_1 \lambda_2^{-1} = 2^{3 m- 2 n} \cdot 5^{n-m} = 1$. 
The unique factorization property of rational numbers then forces 
$3 m - 2 n = n-m = 0$, which yields $m = n = 0$, that is 
condition \eqref{eqn:lt1}. \hfill $\Box$ 
\begin{lemma} \label{lem:lt2} 
Suppose \eqref{eqn:lt1}.   
If $B(\ba; \bp)$ is lower triangular over $\Q(\ba)$ then $p_3 = p_4 = 0$.  
\end{lemma}
{\it Proof}. 
Taking a specialization $\ba = (4,4,4;0,7)$, which is good from 
\eqref{eqn:good}, we have 
\[
B_0 = B_1 = B_2 = 
\begin{pmatrix} 4 & 1 \\ -64/5 & -28/ 5 \end{pmatrix}, \quad 
B_3 = 
\begin{pmatrix} -3/4 & -5/64 \\  1 & 0 \end{pmatrix} \quad 
B_4 = 
\begin{pmatrix} 13/27 & 5/27 \\ -64/27 & -35/27 \end{pmatrix}. 
\]
By formula \eqref{eqn:B(a;p)} we have $B = B_0^{p_0+p_1+p_2} \,  
B_3^{p_3} \, B_4^{p_4} = (B_3 B_4^{-1})^{p_3}$ where condition 
\eqref{eqn:lt1} is also used. 
Observe that the matrix $B_3 B_4^{-1}$ is diagonalized by 
an admissible matrix    
\[
P = 
\begin{pmatrix}
-1/8 & -5/8 \\ 1 & 1 
\end{pmatrix}
\qquad \mbox{as} \qquad  
P^{-1} \, (B_3 \, B_4^{-1}) \, P = \diag \{ \, 2^{-3},  \,\,  - 2^{-3} \cdot 3^3 \, \}, 
\]
and so the lower triangular matrix $B$ is diagonalized as 
$P^{-1} \, B \, P = \diag\{\lambda_1, \, \lambda_2\}$ with 
$\lambda_1 = (2^{-3})^{p_3}$ and $\lambda_2 = (- 2^{-3} \cdot 3^3)^{p_3}$. 
Lemma \ref{lem:eev} then forces $\lambda_1 = \lambda_2$, that is, 
$\lambda_1 \lambda_2^{-1} = (-1)^{p_3} \cdot 3^{-3 p_3} = 1$. 
This immediately shows that $p_3 = -p_4 = 0$.  \hfill $\Box$ 
\begin{lemma} \label{lem:lt3} 
Suppose $p_0+p_1+p_2 = p_3 = p_4 = 0$. 
If $B(\ba; \bp)$ is lower triangular over $\Q(\ba)$ then 
$p_0 = p_1 = p_2 = p_3 = p_4 = 0$, namely, $\bp = \0$.   
\end{lemma}
{\it Proof}. 
Taking a specialization $\ba = (-6,1,1;0,0)$, which is good from 
\eqref{eqn:good}, we have 
\[
B_0 = 
\begin{pmatrix} -6 & 1 \\ -3/2 & -35/ 4 \end{pmatrix}, \quad 
B_1 = B_2 =  
\begin{pmatrix} 1 & 1 \\  -3/2 & -7/4 \end{pmatrix} \quad 
B_3 = B_4 =  
\begin{pmatrix} -11/6 & -2/3 \\ 1 & 0 \end{pmatrix}. 
\]
By formula \eqref{eqn:B(a;p)} we have $B = B_0^{p_0} \, 
B_1^{p_1+p_2} \, (B_3 \, B_4)^{0} = (B_0 B_1^{-1})^{p_0}$ 
where condition $p_0+p_1+p_2 = p_3 = p_4 = 0$ is also used. 
Observe that $B_0 B_1^{-1}$ is diagonalized by an admissible matrix    
\[
P = 
\begin{pmatrix}
-1/2 & -4/3 \\ 1 & 1 
\end{pmatrix}
\qquad \mbox{as} \qquad  
P^{-1} \, (B_0 \, B_1^{-1}) \, P = \diag \{ \, 2^3,  \,\,  - 3^3 \, \}, 
\]
and so the lower triangular matrix $B$ is diagonalized as 
$P^{-1} \, B \, P = \diag\{\lambda_1, \, \lambda_2\}$ with 
$\lambda_1 = 2^{3 p_0}$ and $\lambda_2 = (- 3^3)^{p_0}$. 
Lemma \ref{lem:eev} then forces $\lambda_1 = \lambda_2$, that is, 
$\lambda_1 \lambda_2^{-1} = (-1)^{p_0} \cdot 2^{3 p_0} \cdot 3^{-3 p_0} = 1$. 
This immediately yields $p_0 = -p_1-p_2 = 0$. 
In a similar manner $p_1 = -p_2 - p_0 = 0$ follows from the good 
specialization $\ba = (1,-6,1;0,0)$. 
Thus we have $p_0 = p_1 = p_2 = p_3 = p_4 = 0$.   
\hfill $\Box$ \par\medskip
Lemma \ref{lem:lt} and hence Proposition \ref{prop:lt} follows 
immediately from Lemmas \ref{lem:lt1}, \ref{lem:lt2} and \ref{lem:lt3}. 
Proposition \ref{prop:lt} can be strengthened if the shift vectors $\bp$ are 
restricted somewhat. 
\begin{proposition} \label{prop:lt2} 
Let $\bp = (p_0,p_1,p_2;p_3,p_4) \in \Z_{\ge 0}^3 \times \Z^2$ and $\bc \in \Q^5$. 
Then the specialization 
\[
\hat{r}(a_3,a_4; \bp) := 
r(\ba +\bc; \bp) \big|_{a_0=a_1=a_2=0} \in \Q(a_3, a_4), \qquad 
\ba = (a_0,a_1,a_2;a_3,a_4),  
\]
is well defined.  
If $\hat{r}(a_3,a_4; \bp)$ vanishes in $\Q(a_3, a_4)$ then $\bp$ must be $\0$. 
\end{proposition} 
{\it Proof}. Since $p_0$, $p_1$, $p_2 \ge 0$, one can take a lattice 
path from $\0$ to $\bp$ in such a manner that if $i_k \in \{0, 1, 2\}$ then 
$\ve_k = +$ in \eqref{eqn:pm}, so that the product \eqref{eqn:A(a;p)} contains 
no contiguous matrices of type \eqref{eqn:cm3} in Table \ref{tab:cm}.  
In view of the denominators of contiguous matrices of 
types \eqref{eqn:cm1}, \eqref{eqn:cm2}, \eqref{eqn:cm4}, any irreducible 
factor of the denominator 
of $A(\ba + \bc; \bp)$ must be either $s(\ba) + \mbox{a rational number}$,  
or $a_j - a_i + \mbox{a rational number}$, for some $i = 0,1,2$, $j = 3, 4$. 
Thus one can safely take the specialization $\hat{A}(a_3, a_4; \bp) := 
A(\ba + \bc; \bp) \big|_{a_0=a_1=a_2= 0}$ as a matrix over $\Q(a_3,a_4)$, 
the principal part of which is independent of $\bc$ and given by 
$\hat{B}(a_3, a_4; \bp) = \hat{B}_0^{p_0} \hat{B}_1^{p_1}  \hat{B}_2^{p_2} 
\hat{B}_3^{p_3} \hat{B}_4^{p_4}$, where    
\[
\hat{B}_i  := 
\begin{pmatrix} 
0 & 1 \\[2mm] 
0 & - a_3 a_4 (a_3+a_4)^{-1} 
\end{pmatrix},  
\quad  i = 0, 1, 2; \qquad 
\hat{B}_i  :=  
\begin{pmatrix} 
a_i^{-1} & -a_i^{-3}(a_3+a_4) \\[2mm]  
0 & a_i^{-2}(a_3+ a_4) 
\end{pmatrix},  
\quad i = 3, 4,  
\]
are derived from formulas \eqref{eqn:pp} in Table \ref{tab:pp} by 
putting $a_0=a_1=a_2=0$. 
Since $\hat{B}_0 = \hat{B}_1 = \hat{B}_2$, one has 
$\hat{B}(a_3, a_4; \bp) = \hat{B}_0^p \hat{B}_3^{p_3} \hat{B}_4^{p_4}$ 
with $p := p_0+p_1+p_2$. 
In place of \eqref{eqn:good} the goodness condition is now 
$a_3 a_4(a_3+a_4) \neq 0$.  
Take a good specialization $a_3 = -2$, $a_4 = 3$. 
After some manipulations we have 
\begin{alignat*}{2}
\hat{B}(-2, 3; \bp) &= 
\begin{pmatrix}
\, (-2)^{-p_3} \cdot 3^{-p_4}  \,\, & 2^{-2 p_3-1} \cdot 3^{-2 p_4-1} 
\{1-(-2)^{p_3} \cdot 3^{p_4} \} \, \\[1mm]
\, 0  \,\, & 2^{-2 p_3} \cdot 3^{-2 p_4} \,  
\end{pmatrix}, 
\qquad & p &= 0, \\[2mm] 
\hat{B}(-2, 3; \bp) &= 
\begin{pmatrix} 
\, 0 \,\, & 2^{p-2 p_3-1} \cdot 3^{p-2 p_4-1} \, \\[1mm] 
\, 0 \,\, & 2^{p-2 p_3} \cdot 3^{p-2p_4} \,  
\end{pmatrix}, 
\qquad & p &\ge 1.   
\end{alignat*}
If $\hat{r}(a_3, a_4; \bp) = 0$ in $\Q(a_3,a_4)$ then 
$\hat{A}(a_3, a_4; \bp)$ and so $\hat{B}(-2, 3; \bp)$ must be 
lower triangular. 
This forces $p=p_0+p_1+p_2=0$ and 
$(-2)^{p_3} \cdot 3^{p_4} = 1$. 
Since $p_0$, $p_1$, $p_2 \ge 0$, we must have 
$p_0 = p_1 = p_2 = p_3 = p_4 = 0$ and hence the proposition is proved.  
\hfill $\Box$ \par\medskip
Proposition \ref{prop:lt2} is used to define specializations of ${}_3F_2(1)$ continued 
fractions in \cite{EI}. 
\section{Linear Independence} \label{sec:lin-indep} 
The aim of this section is to establish the linear independence of 
$\hgh(\ba+\bp)$ and $\hgh(\ba+\bq)$ for any distinct $\bp$, $\bq \in \Z^5$ 
(Proposition \ref{prop:indep-pq}); this is indispensable for the existence 
and uniqueness of solutions to the linear 
system \eqref{eqn:one-p}-\eqref{eqn:one-q}.  
Once this is done then the proof of our main theorem (Theorem \ref{thm:3tr}) 
is an easy task and hence completed in this and next sections.    
\par
We begin by dealing with an important special case where $\bp = \0$ 
and $\bq = \1$.       
For this case three-term relation \eqref{eqn:one-two} can be used to 
show the following fundamental lemma. 
\begin{lemma} \label{lem:one}
$\hgh(\ba)$ and $\hgh(\ba + \1)$ are linearly independent over $\C(\ba)$. 
\end{lemma}
{\it Proof}. 
To get a contradiction, suppose that $\hgh(\ba)$ and 
$\hgh(\ba + \1)$ are linearly dependent over $\C(\ba)$. 
Then there exists a rational function $R(\ba) \in \C(\ba)$ such that 
$\hgh(\ba + \1) = R(\ba) \, \hgh(\ba)$. 
This implies $\hgh(\ba +\2) = R(\ba +\1) \, R(\ba) \, \hgh(\ba)$. 
Substituting these into \eqref{eqn:one-two} yields  
\begin{equation} \label{eqn:one1}
\{s(\ba) -2 \} \, R(\ba +\1) \, R(\ba) + \psi(\ba) \, R(\ba) - 
\varphi_3(\ba) = 0, 
\end{equation}
where $\psi(\ba) = b_1 b_2 - \varphi_2(\ba) - \varphi_1(\ba) -1$. 
Putting $R(\ba) = P(\ba)/Q(\ba)$ with $P(\ba)$, $Q(\ba) \in \C[\ba]$ 
and multiplying \eqref{eqn:one1} by $Q(\ba +\1) \, Q(\ba)$, we have  
a polynomial equation 
\begin{equation} \label{eqn:one2}
\{s(\ba) -2 \} \, P(\ba +\1) \, P(\ba) + 
\psi(\ba) \, P(\ba) \, Q(\ba + \1) - 
\varphi_3(\ba) \, Q(\ba +\1) \, Q(\ba) = 0.  
\end{equation}
\par
If $d_p := \deg P(\ba)$ and $d_q := \deg Q(\ba)$, the three terms 
in the left side of \eqref{eqn:one2} have degrees $2 d_p +1$, 
$d_p + d_q + 2$, $2 d_q + 3$, at least two of which must coincide.  
This implies $d_p = d_q +1$ and hence the three degrees must be equal. 
Thus the top homogeneous component of \eqref{eqn:one2} gives 
\begin{equation} \label{eqn:one3}
s(\ba) \, \bar{P}^2(\ba) + \{ b_1 b_2 - \varphi_2(\ba) \} \, \bar{P}(\ba) 
\, \bar{Q}(\ba)  - \varphi_3(\ba) \, \bar{Q}^2(\ba) = 0,  
\end{equation}
where $\bar{P}(\ba)$ and $\bar{Q}(\ba)$ are the top homogeneous 
components of $P(\ba)$ and $Q(\ba)$ respectively. 
We write $\bar{P}(\ba) = p(\ba) \, r(\ba)$ and 
$\bar{Q}(\ba) = q(\ba) \, r(\ba)$ with $p(\ba)$, $q(\ba)$, $r(\ba) 
\in \C[\ba]$, where $p(\ba)$ and $q(\ba)$ have no factors in common. 
Then \eqref{eqn:one3} and $d_p = d_q + 1$ yield  @ 
\begin{equation} \label{eqn:one4} 
s(\ba) \, p^2(\ba) + \{ b_1 b_2 - \varphi_2(\ba) \} \, p(\ba) 
\, q(\ba)  - \varphi_3(\ba) \, q^2(\ba) = 0, \quad 
\deg p(\ba) = \deg q(\ba) + 1.    
\end{equation}
\par
The former equation in \eqref{eqn:one4} implies 
$q(\ba) \,|\, s(\ba) \, p^2(\ba)$ and 
$p(\ba) \, | \, \varphi_3(\ba) \, q^2(\ba)$. 
Since $p(\ba)$ and $q(\ba)$ have no factors in common, we must have  
$q(\ba) \,|\, s(\ba)$ and $p(\ba) \,|\, \varphi_3(\ba) = a_1 a_2 a_3$. 
These division relations and the latter equation in \eqref{eqn:one4} 
lead to a dichotomy  
\begin{alignat*}{3} 
\mbox{(i)} & \qquad & q(\ba) &= c_q, \qquad & p(\ba) &= c_p \, a_i,  \\
\mbox{(ii)} & \qquad & q(\ba) &= c_q \, s(\ba), \qquad & p(\ba) &= c_p \, a_i \, a_j, 
\end{alignat*}    
where $c_p$, $c_q \in \C^{\times}$ and $\{i, j, k\} = \{1,2,3\}$. 
The former equation in \eqref{eqn:one4} then gives   
\begin{alignat*}{2} 
\mbox{(i)} & \qquad & 
c_p^2 \, a_i \, s(\ba) + c_p c_q \, \{b_1 b_2 - \varphi_2(\ba) \} - 
c_q^2 \, a_j \, a_k &= 0, \\
\mbox{(ii)} & \qquad & 
c_p^2 \, a_i \, a_j + c_p c_q \, \{b_1 b_2 - \varphi_2(\ba) \} - 
c_q^2 \, a_k \, s(\ba) &= 0, 
\end{alignat*}    
neither of which is feasible as an equation in $\C[\ba]$. 
This contradiction proves the lemma. \hfill $\Box$  
\begin{remark} \label{rem:one} 
We make two remarks about Lemma \ref{lem:one}. 
\begin{enumerate} 
\item Vidunas \cite{Vidunas} showed the linear independence over 
$\C(\ba;z)$ of ${}_2F_1(\ba;z)$ and ${}_2F_1(\ba+\be_1;z)$ with 
$\be_1 := (1,0;0)$, using Kummer's ${}_2F_1(-1)$ formula as in 
\cite[(14)]{Vidunas}.  
Lemma \ref{lem:one} might be proved in a similar spirit, 
while our proof here is purely algebraic and totally independent of 
special-value formulas, which we believe is of its own merit and interest;   
e.g., availability of field extensions as in item (2) below and potential 
applicability to restricted $\hgF(\ba)$'s.            
\item From the way in which Lemma \ref{lem:one} is proved it is evident that 
the lemma remains true over $K(\ba)$ for any function field 
$K \supset \Q$ in $\ba$ satisfying condition \eqref{eqn:transc}. 
\end{enumerate} 
\end{remark}
\par
Proposition \ref{prop:lt}, Lemmas \ref{lem:solve} and \ref{lem:one} are put 
together to establish the following. 
\begin{proposition} \label{prop:indep-pq} 
If $\bp$, $\bq \in \Z^5$ are distinct integer vectors, then  
\begin{enumerate}
\item the determinant $\vD(\ba; \bp, \bq)$ in \eqref{eqn:nonzero} is 
non-vanishing in $\Q(\ba);$  
\item $\hgh( \ba + \bp )$ and $\hgh( \ba + \bq)$ are linearly 
independent over $K(\ba)$ for any function field $K \supset \Q$ in 
$\ba$ satisfying condition \eqref{eqn:transc}$;$ 
\item there exist unique rational functions $u(\ba)$, $v(\ba) \in \Q(\ba)$ such that 
\begin{equation} \label{eqn:3tr2}
\hgh(\ba) = u(\ba) \, \hgh(\ba+\bp) + v(\ba) \, \hgh(\ba+\bq),   
\end{equation}
where $u(\ba)$ and $v(\ba)$ are common for all choices of $\hgh(\ba)$ and 
explicitly given by   
\begin{subequations} \label{eqn:u(a;p,q)} 
\begin{align}
u(\ba) &= \phantom{-} \dfrac{r(\ba; \bq)}{\vD(\ba;\bp,\bq)} 
= \phantom{-} \dfrac{r(\ba; \bq)}{\delta(\ba; \bp) \cdot r(\ba+\bp; \bq-\bp)}, 
\label{eqn:u(a;p,q)u} \\[2mm] 
v(\ba) &= -\dfrac{r(\ba; \bp)}{\vD(\ba; \bp, \bq)} 
= -\dfrac{r(\ba; \bp)}{\delta(\ba; \bp) \cdot r(\ba+\bp; \bq-\bp)}, 
\label{eqn:u(a;p,q)v}
\end{align} 
\end{subequations}
where $\delta(\ba; \bp)$ is defined in \eqref{eqn:det}.  
\end{enumerate} 
\end{proposition}
{\it Proof}. 
Proposition \ref{prop:lt} implies that $r(\ba + \bp; \bq-\bp)$ is nonzero 
in $\Q(\ba)$, then assertion (1) follows from equation \eqref{eqn:solve} in 
Lemma \ref{lem:solve}. 
By assertion (1) the matrix in \eqref{eqn:base-change} is invertible over 
$\Q(\ba)$ and hence over $K(\ba)$, so assertion (2) follows from   
Lemma \ref{lem:one} together with item (2) of Remark \ref{rem:one}. 
The second equations in \eqref{eqn:u(a;p,q)u} and \eqref{eqn:u(a;p,q)v} 
come from Lemma \ref{lem:solve}.  
Converting equation \eqref{eqn:base-change} one has three-term relation 
\eqref{eqn:3tr2} with $u(\ba)$ and $v(\ba)$ given by \eqref{eqn:u(a;p,q)}. 
The uniqueness of $u(\ba)$ and $v(\ba)$ in \eqref{eqn:3tr2} is an easy 
consequence of assertion (2). 
Assertion (3) is thus established. \hfill $\Box$ \par\medskip
Summarizing the discussions in \S\ref{sec:cm} and in this section 
we have the following recipe to determine the coefficients of three-term 
relation \eqref{eqn:3tr}. 
It works effectively on computers.   
\begin{recipe} \label{recipe}
Given any distinct $\bp$, $\bq \in \Z^5$, take a lattice path from 
$\0$ to $\bp$ and calculate the successive product \eqref{eqn:A(a;p)} of 
contiguous matrices along it to produce the connection matrix $A(\ba;\bp)$. 
From it extract its upper components $r_1(\ba; \bp)$ and $r(\ba; \bp)$ 
as in \eqref{eqn:ur}.  
With $\bq$ in place of $\bp$, proceed exactly in the same manner 
to get $r_1(\ba; \bq)$ and $r(\ba; \bq)$. 
Then the coefficients $u(\ba)$ and $v(\ba)$ of three-term relation 
\eqref{eqn:3tr} are given by the first equations in \eqref{eqn:u(a;p,q)u} and 
\eqref{eqn:u(a;p,q)v}. 
Alternatively one may use the second equations in \eqref{eqn:u(a;p,q)u} and 
\eqref{eqn:u(a;p,q)v}, in which case one should also take the product 
along a lattice path from $\0$ to $\bq-\bp$ to deduce $r(\ba+\bp; \bq-\bp)$. 
\end{recipe}
\par
With Proposition \ref{prop:indep-pq} and Recipe \ref{recipe} 
all assertions of Theorem \ref{thm:3tr} have been established except for 
assertion (5), which is treated in the next section 
(see Proposition \ref{prop:symmetry}).     
\section{Symmetry} \label{sec:symmetry} 
The uniqueness and simultaneousness for the three-term relations provide  
themselves with a group symmetry of order seventy-two.  
To see this let 
\[
G := \langle \sigma_0, \, \sigma_1, \, \sigma_2, \, \sigma_3, \, \sigma_4 \rangle 
\qquad \mbox{with} \quad \sigma_3 :=  \tau_1, \quad \sigma_4 := \tau_2
\]
be the group of affine transformations generated by the involutions in 
Table \ref{tab:inv}. 
The Saalsch\"{u}tzian index $s(\ba)$ is preserved by the generators and 
hence by the group $G$. 
Note also that   
\[
\bc := (2/3, 2/3, 2/3; 1, 1)
\]
is the unique fixed point of the generators and hence of $G$ with null 
Saalsch\"{u}tzian index. 
\par
Taking the linear parts of affine transformations induces 
a surjective group homomorphism  
\begin{equation*} 
G \rightarrow \bar{G} := 
\langle \bar{\sigma}_0, \, \bar{\sigma}_1, \, \bar{\sigma}_2, \, 
\bar{\sigma}_3, \, \bar{\sigma}_4 \rangle, \qquad \sigma \mapsto \bar{\sigma},  
\end{equation*}
where $\bar{\sigma}$ stands for the linear part of an affine $\sigma \in G$.  
This is just an isomorphism because  
\[
\bar{\sigma} = t^{-1} \cdot \sigma \cdot t \qquad 
\mbox{with} \quad t : \ba \mapsto \ba + \bc \quad 
\mbox{being the parallel translation by $\bc$}.    
\]
Each $\bar{\sigma}_i$ is a $\Z$-linear transformation of determinant 
$\det \bar\sigma_i = \pm 1$, actually $+1$ for $i = 0,1,2$ and $-1$ for 
$i = 3,4$, so that every $\bar{\sigma} \in \bar{G}$ 
maps the lattice $\Z^5$ isomorphically onto itself. 
\par
Put $u(\ba; \bp, \bq) := r(\ba; \bq)/\vD(\ba; \bp, \bq)$. 
Note that $u(\ba) = u(\ba; \bp, \bq)$ and $v(\ba) = u(\ba; \bq, \bp)$ 
in formulas \eqref{eqn:u(a;p,q)}.   
We show that $u(\ba; \bp, \bq)$ enjoys the following $G$-covariance.  
\begin{lemma} \label{lem:g-action} 
For any distinct $\bp$, $\bq \in \Z^5$ and any $\sigma \in G$ we have   
\begin{equation} \label{eqn:g-action} 
u(\ba; \bar{\sigma}(\bp), \bar{\sigma}(\bq)) = u(\sigma^{-1}(\ba); \bp, \bq). 
\end{equation}
\end{lemma}
{\it Proof}. 
In three-term relation \eqref{eqn:3tr2} we take $\hgh(\ba) = 
\hgf(\sigma_i(\ba))$ for each $i = 0,1,2,3,4$. 
Formula \eqref{eqn:3tr2} with $\hgh$, $\ba$, $\bp$ and $\bq$ 
replaced by $\hgf$, $\sigma_i(\ba)$, $\bar{\sigma}_i(\bp)$ and  
$\bar{\sigma}_i(\bq)$ yields  
\begin{align*}
\hgh(\ba) 
&= \hgf(\sigma_i(\ba)) \\ 
&= u(\sigma_i(\ba); \bar{\sigma}_i(\bp), \bar{\sigma}_i(\bq)) \, 
\hgf(\sigma_i(\ba) + \bar{\sigma}_i(\bp)) + 
u(\sigma_i(\ba); \bar{\sigma}_i(\bq), \bar{\sigma}_i(\bp)) \, 
\hgf(\sigma_i(\ba) + \bar{\sigma}_i(\bq)) \\ 
&= u(\sigma_i(\ba); \bar{\sigma}_i(\bp), \bar{\sigma}_i(\bq)) \, 
\hgf(\sigma_i(\ba + \bp)) + 
u(\sigma_i(\ba); \bar{\sigma}_i(\bq), \bar{\sigma}_i(\bp)) \, 
\hgf(\sigma_i(\ba + \bq)) \\
&= u(\sigma_i(\ba); \bar{\sigma}_i(\bp), \bar{\sigma}_i(\bq)) \, 
\hgh(\ba + \bp) + 
u(\sigma_i(\ba); \bar{\sigma}_i(\bq), \bar{\sigma}_i(\bp)) \, 
\hgh(\ba + \bq),   
\end{align*}
which must coincide with \eqref{eqn:3tr2} by the uniqueness of 
three-term relation. 
Thus we have 
\[ 
u(\sigma_i(\ba); \bar{\sigma}_i(\bp), \bar{\sigma}_i(\bq)) = 
u(\ba; \bp, \bq), \quad  
u(\sigma_i(\ba); \bar{\sigma}_i(\bq), \bar{\sigma}_i(\bp)) = 
u(\ba; \bq, \bp), \quad i = 0,1,2,3,4,   
\] 
and hence $u(\sigma(\ba); \bar{\sigma}(\bp), \bar{\sigma}(\bq)) 
= u(\ba; \bp, \bq)$ for every $\sigma \in G$. 
Formula \eqref{eqn:g-action} then follows from the last equation 
by replacing $\ba$ with $\sigma^{-1}(\ba)$. 
\hfill $\Box$ \par\medskip 
We are interested in the group structure of $G$ or $\bar{G}$.  
It is easy to see that     
\[
\rho_i := \tau_i \, \sigma_0 \, \tau_i \, \sigma_0 \, \tau_i = 
\tau_i \, \sigma_i \, \tau_i \, \sigma_i \, \tau_i \in G, \qquad i = 1,2,   
\]
are involutions such that $\sigma_i = \rho_i \, \sigma_0 \, \rho_i$ 
and $\tau_i = \sigma_0 \, \rho_i \, \sigma_0$ for $i = 1, 2$, 
so that we have  
\[
G = \langle \sigma_0, \, \rho_1, \, \rho_2, \, \tau_1, \, \tau_2 
\rangle = \langle \sigma_0, \, \rho_1, \, \rho_2 \rangle. 
\]    
\begin{lemma} \label{lem:gs} 
The group structure of $\bar{G}$ is given by    
\begin{equation} \label{eqn:gs}
\bar{G} = \langle \bar{\sigma}_0 \rangle \ltimes 
(\langle \bar{\rho}_1, \, \bar{\rho}_2 \rangle \times 
\langle \bar{\tau}_1, \, \bar{\tau}_2 \rangle) 
\cong S_2 \ltimes (S_3 \times S_3),  
\end{equation}
in particular $\bar{G}$ is a group of order $2! \times (3! \times 3!) = 72$. 
Viewed as acting on the real vector space $\R^5$ with coordinates 
$\bp = (p_0,p_1,p_2; q_1, q_2)$, the group $\bar{G}$ has a fundamental domain 
\begin{equation} \label{eqn:fd} 
p_0 \ge p_1 \ge p_2, \quad q_2 \ge q_1 \ge p_0 - p_1. 
\end{equation}
\end{lemma} 
{\it Proof}.  
As a linear basis of $\R^5$ we take five vectors   
\begin{alignat*}{3}  
\bu_1 &:= (1/3, -2/3, 1/3; 0, 0), \quad & \bv_1 &:= (1/3, 1/3, 1/3; 1, 0), \quad 
  \bw &:= (-1/3,-1/3,-1/3; 0, 0),               \nonumber                     \\
\bu_2 &:= (1/3, 1/3, -2/3; 0, 0), \quad & \bv_2 &:= (1/3, 1/3, 1/3; 0, 1), \quad 
      &                                 &        
\end{alignat*}
along with two auxiliary vectors such that $\bu_0+\bu_1+\bu_2 = \0$ and 
$\bv_0+\bv_1+\bv_2 = \0$, namely,    
\[
\bu_0 := (-2/3, 1/3, 1/3; 0, 0), \qquad \bv_0 := (-2/3, -2/3, -2/3; -1, -1). 
\]
The linear action of $\bar{G}$ on $\R^5$ is faithfully represented 
by permutations on these vectors: 
\[
\bar{\rho}_i \, : \, \bu_0 \leftrightarrow \bu_i, \quad  
\bar{\tau}_i \, : \, \bv_0 \leftrightarrow \bv_i, \quad i = 1, 2; \qquad     
\bar{\sigma}_0 \, : \, \bu_0 \leftrightarrow \bv_0, \,\,  
\bu_1 \leftrightarrow \bv_1, \,\, \bu_2 \leftrightarrow \bv_2,  
\]
where the vectors left invariant are not indicated. 
Thus we have two actions    
\[ 
S_3 \cong \langle \bar{\rho}_1, \, \bar{\rho}_2 \rangle \curvearrowright 
U := \R \, \bu_1 \oplus \R \, \bu_2, \qquad 
S_3 \cong \langle \bar{\tau}_1, \, \bar{\tau}_2 \rangle \curvearrowright 
V:= \R \, \bv_1 \oplus \R \, \bv_2,    
\]
which are commutative and permuted by $\bar{\sigma}_0$, together with 
a line $\R \, \bw$ fixed pointwise by $\bar{G}$. 
These observations clearly show that $\bar{G}$ has the group structure 
\eqref{eqn:gs}.    
If $\bp$ is represented as 
\[
\bp = x_1 \, \bu_1 + x_2 \, \bu_2 + y_1 \, \bv_1 + y_2 \, \bv_2 + z \bw \in 
U \oplus V \oplus \R \, \bw, 
\]
then $x_1 = p_0 - p_1$, $x_2 = p_0 - p_2$, $y_1 = q_1$, $y_2 = q_2$ and 
$z = s(\bp)$. 
As a fundamental domain of 
$\langle \bar{\rho}_1, \, \bar{\rho}_2 \rangle \times 
\langle \bar{\tau}_1, \, \bar{\tau}_2 \rangle$ we can take 
$D_0 := \{\, x_2 \ge x_1 \ge 0, \,\, y_2 \ge y_1 \ge 0 \,\}$. 
Note that $\bar{\sigma}_0$ maps $D_0$ onto itself while swapping  
two conditions $x_2 \ge x_1 \ge 0$ and $ y_2 \ge y_1 \ge 0$. 
Thus we can impose condition $y_1 \ge x_1$ to get a fundamental domain 
$D :=  \{\, x_2 \ge x_1 \ge 0, \,\, y_2 \ge y_1 \ge 0, \, \, 
y_1 \ge x_1 \,\}$ of the whole group $\bar{G}$. 
In terms of the original coordinates of $\bp$ the domain $D$ 
is given by \eqref{eqn:fd}.  \hfill $\Box$ \par\medskip
Putting Lemmas \ref{lem:g-action} and \ref{lem:gs} together we have 
the following. 
\begin{proposition} \label{prop:symmetry} 
Three-term relation \eqref{eqn:3tr} admits a 
$S_2 \ltimes (S_3 \times S_3)$-symmetry 
\begin{equation} \label{eqn:symmetry}
\hgh(\ba) = {}^{\sigma}\!u(\ba) \, \hgh(\ba+\bar{\sigma}(\bp)) + 
{}^{\sigma}\!v(\ba) \, \hgh(\ba+\bar{\sigma}(\bq)), \qquad 
\sigma \in G, 
\end{equation} 
where ${}^{\sigma}\!w(\ba) := w(\sigma^{-1}(\ba))$ is the induced 
action of $\sigma$ on a function $w(\ba)$.   
\end{proposition} 
\begin{remark} \label{rem:symmetry} 
A couple of remarks are in order at this stage. 
\begin{enumerate}
\item Formula \eqref{eqn:symmetry} is of particular interest when $\bq$ is 
$\be :=(1,1,1;0,0)$ or a nonzero integer multiple of it, because in that 
case $\bq$ is $\bar{G}$-invariant and hence \eqref{eqn:symmetry} becomes
\[
\hgh(\ba) = {}^{\sigma}\!u(\ba) \, \hgh(\ba+\bar{\sigma}(\bp)) + 
{}^{\sigma}\!v(\ba) \, \hgh(\ba+\bq), \qquad 
\sigma \in G,  
\]
so without loss of generality $\bp$ may lie within a fundamental domain 
of $\bar{G}$ as in \eqref{eqn:fd}. 
\item In \S \ref{sec:cm} the contiguous and connection matrices 
$A_i^{\ve}(\ba)$ and $A(\ba; \bp)$ were formulated in terms of the basis 
$\bhgh(\ba) = {}^t( \hgh(\ba), \, \hgh(\ba+\1) )$ as in 
\eqref{eqn:vec-h}, where $\1$ was not $\bar{G}$-invariant.   
From the viewpoint of symmetry it is better to reformulate them 
in term of an alternative basis $\tilde{\bhgh}(\ba) := 
{}^t( \hgh(\ba), \, \hgh(\ba+\be) )$. 
The revised contiguous and connection matrices are   
\[
\tilde{A}_i^{\ve}(\ba) = 
P(\ba + \ve \, \be_i) \, A_i^{\ve}(\ba) \, P(\ba)^{-1}, \qquad 
\tilde{A}(\ba;\bp) = P(\ba + \bp) \, A(\ba; \bp) \, P(\ba)^{-1},   
\]
where $P(\ba)$ is the invertible matrix over $\Q(\ba)$ such that  
$\tilde{\bhgh}(\ba) = P(\ba) \, \bhgh(\ba)$. 
An advantage of this base change is that $\tilde{A}(\ba;\bp)$ gains  
the nice $G$-covariance property    
\[
\tilde{A}(\ba; \bar{\sigma}(\bp)) = \tilde{A}(\sigma^{-1}(\ba); \bp), 
\qquad \sigma \in G,      
\]
but unfortunately explicit formula for $\tilde{A}_i^{\ve}(\ba)$ 
is too complicated to be presented here. 
\item Let $\rho_3 := \tau_1 \tau_2 \tau_1 = \tau_2 \tau_1 \tau_2 \in G$. 
Then $G$ contains a subgroup $G_0 := \langle \rho_1, \, \rho_2 \rangle 
\times \langle \rho_3 \rangle \cong S_3 \times S_2$ that acts on 
$\hgf(\ba)$ trivially permuting its numerator or denominator 
parameters. 
The left quotient set $G_0 \backslash G$ is then represented by 
$\sigma_i$, $i = 0,1,2,3,4$, and the unit element, so the $G$-orbit 
of $\hgf(\ba)$ is exactly the six functions in 
\eqref{eqn:companion} up to the trivial symmetry $G_0$.     
\end{enumerate} 
\end{remark} 
\section{Concluding Remarks} \label{sec:concluding}
Theorem \ref{thm:3tr} contains five key words:       
{\sl linear independence}, {\sl existence}, {\sl uniqueness}, 
{\sl simultaneousness} and {\sl symmetry} for general three-term 
relations with arbitrary integer shifts.  
Among them existence and uniqueness are easy consequences of 
linear independence (on the basis of three-term contiguous 
relations in the narrow sense), while symmetry is a direct 
corollary to uniqueness and simultaneousness. 
In this sense linear independence and simultaneousness are 
the most fundamental concepts in our discussions. 
\par
Another simple corollary to the linear independence is the 
{\sl impossibility} of representing unrestricted $\hgF(1)$ series as 
products of gamma functions like  
\[
\hgF(\ba; 1) = C \cdot d^{K(\sba)} 
\dfrac{\prod_{i=1}^I \vG(L_i(\ba))}{\prod_{j=1}^J \vG(M_j(\ba))}, 
\]
where $C$, $d \in \C^{\times}$ are nonzero constants, 
$K(\ba)$, $L_i(\ba)$, $M_j(\ba)$ are affine polynomials over $\C$ 
in $\ba$ such that the homogeneous linear parts of $L_i(\ba)$ 
and $M_j(\ba)$ have coefficients in $\Q$. 
Indeed, if such an identity existed then Proposition \ref{prop:indep-pq} 
and the recursion formula for the gamma function $\vG(z+1) = z \, \vG(z)$ 
would readily lead to a contradiction. 
This result remains true over certain function fields $K \supset \Q$ in 
place of $\C$.  
But an even stronger result is already known by Wimp 
\cite[Theorems 3 and 4]{Wimp} based on the concept of {\sl unicial} functions. 
The depth of his result is due to the use of a transcendental method, 
that is, asymptotic analysis of linear recurrence equations. 
By contrast, our method relies only on linear algebra and hence 
is more elementary.     
\par
We can understand and simplify Wimp's difficult arguments in 
\cite[\S\S 2-3]{Wimp} from the general framework of this article.  
His recurrence equation \cite[equation (13)]{Wimp} just stems from 
the three-term relation \eqref{eqn:3tr} with shift vectors $\bp = (1,1,1;2,1)$ 
and $\bq = - \bp$, which can easily be calculated by applying Recipe 
\ref{recipe} to these vectors.  
It is then trivial that his sequence $C(n)$ satisfies 
\cite[equation (13)]{Wimp}. 
His statements that the related sequences in 
\cite[formulas (15), (16), (17)]{Wimp} also satisfy the same recurrence 
are consequences of simultaneousness and symmetry. 
Indeed, the statements for the sequences $C_4(n)$ and $C_5(n)$ come from the 
symmetries $\tau_1$ and $\tau_2$ respectively, while those for $C_h(n)$, 
$h = 1,2,3$, are due to $\sigma_0$ followed by permutations of numerator 
or denominator parameters. 
\par
Lastly but most importantly, three-term contiguous relations give rise to 
three-term recurrence relations, which in turn generate continued fractions, 
so the results and methods of this article are fundamental in developing 
a general theory of $\hgF(1)$ continued fractions in \cite{EI}, which 
establishes an infinite number of continued fraction expansions having   
exact error term estimates and enjoying rapid convergences  
for ratios of ${}_3F_2(1)$ hypergeometric series.   
\par\vspace{3mm} \noindent
{\bf Acknowledgment}. 
This work is supported by Grant-in-Aid for Scientific Research, 
JSPS, 16K05165 (C).                  
\bibliographystyle{plain}
\bibliography{3tr3f2ref}
\end{document}